\documentclass[11pt]{article}  
\usepackage{amsmath}
\usepackage{graphicx} 
\usepackage{amssymb} 
\usepackage{amsthm}
\usepackage{subcaption} 
\usepackage{tikz}
\usepackage{authblk}

\usepackage{comment}
\usepackage{xcolor}

\usepackage{hyperref}
\hypersetup{
    colorlinks,
    linkcolor={red!50!black},
    citecolor={blue!50!black},
    urlcolor={blue!80!black}
}

\usepackage[margin=1.5in]{geometry}
\usetikzlibrary{shapes}
\usetikzlibrary{positioning}
\usetikzlibrary{arrows}
\usetikzlibrary{decorations.markings}
\usetikzlibrary{arrows.meta,decorations.markings}

\theoremstyle{plain}

\newtheorem{theorem}{Theorem}[section]
\newtheorem{lemma}[theorem]{Lemma}
\newtheorem{corollary}[theorem]{Corollary}

\DeclareMathOperator{\tr}{tr} 

\newcommand{\ssym}[0]{spectrally symmetric }
\newcommand{\Hs}[0]{H_{\frac{\pi}{3}}}
\newcommand{\ssymend}[0]{spectrally symmetric}
\newcommand{\W}[0]{\mathcal{W}}
\newcommand{\T}[0]{\mathcal{T}}
\renewcommand{\P}[0]{\mathcal{P}}
\renewcommand{\S}[0]{\mathcal{S}}

\title{Spectrally symmetric orientations of graphs}

\author[1]{Saieed Akbari\thanks{The research visit of S. Akbari at Simon Fraser University was supported in part by the ERC Synergy grant (European Union, ERC, KARST, project number 101071836).}}
\author[2]{Jonathan Aloni\thanks{Supported by the NSERC Undergraduate Student Research Award.}} 
\author[3]{Maxwell Levit\thanks{Corresponding Author: Maxwell.Levit@gmail.com Postdoctoral fellowship at SFU supported through NSERC Discovery Grants R832714 and R611368 (Canada), and the ERC Synergy grant KARST (European Union, ERC, KARST, project number 101071836).} }
\author[2]{Bojan Mohar\thanks{Supported in part by the NSERC Discovery Grant R832714 (Canada), and in part by the ERC Synergy grant KARST (European Union, ERC, KARST, project number 101071836).
On leave from FMF, Department of Mathematics, University of Ljubljana.}}
\author[2]{Steven Xia}
\affil[1]{Sharif University of Technology}
\affil[2]{Simon Fraser University}
\affil[3]{York University}

\date{}

\begin{document}
\maketitle
\begin{abstract}
A graph is spectrally symmetric if its adjacency eigenvalues are symmetric about zero. A spectrally symmetric graph must be bipartite and hence, by Mantel's theorem, it contains at most $n^2/4$ edges. Using a natural measure of density, this says that all spectrally symmetric graphs have density at most $\tfrac{1}{2}$. When it comes to digraphs, the corresponding results are much less predictable. 
Here spectral symmetry is defined using the eigenvalues of the Hermitian adjacency matrices of the second kind, a natural choice for the spectral theory of digraphs. In this oriented graph setting, spectrally symmetric orientations need not be bipartite. Indeed we provide examples with density as large as $\tfrac{31}{36}$. Nevertheless, in our main result we are able to enforce an upper bound of $\tfrac{10}{11}$ on the density of any spectrally symmetric oriented graph. This bound can be viewed as a spectral analog of Mantel's theorem for digraphs. We then give another necessary conditions for spectral symmetry, and provide a general recipe for constructing infinite families of examples using 1-sums.
\end{abstract}

\section{Introduction}

In this paper all graphs and digraphs are loopless. Let $G$ be a graph, and let $A=A(G)$ be the adjacency matrix of $G$. It is well known that $G$ is bipartite if and only if the spectrum of $A$ is symmetric about 0. 
Now, let $D$ be a digraph in which $a(i,j)$ denotes the number of arcs from $i$ to $j$. For $\theta\in (0,\pi]$ the \emph{Hermitian adjacency matrix}, $H_\theta(D)$ is defined by \[H_\theta(D)_{i,j}=e^{i\theta} a(i,j)+e^{-i\theta} a(j,i).\]

\begin{figure} 
    \centering
    \includegraphics[scale=0.20, angle=270]{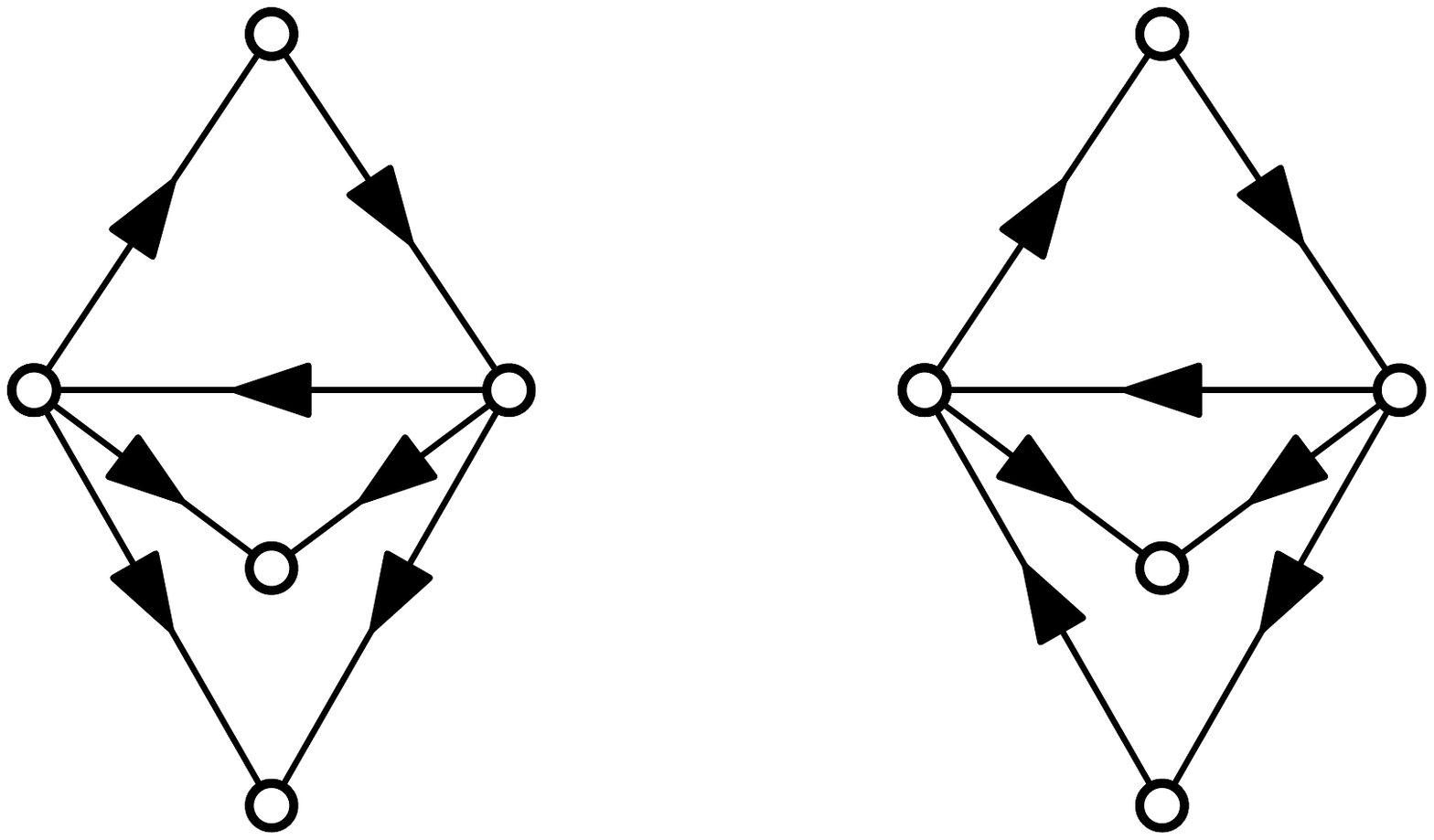}
    \caption{Two \ssym orientations of $K_5-E(K_3)$.}
    \label{fig_k23}
\end{figure}

A digraph is \emph{bipartite} if its underlying graph is bipartite. If $D$ is a bipartite digraph and $\lambda$ an eigenvalue of $H_\theta(D)$, one can construct an eigenvector of $H_\theta(D)$ with eigenvalue $-\lambda$ in exactly the same way one does for the adjacency matrix of an undirected graph. So all bipartite digraphs have symmetric $H_\theta$-spectrum.  The converse is false: Figure \ref{fig_k23} shows an orientation of $K_5-E(K_3)$ which was shown to have symmetric $H_{\frac{\pi}{3}}$-spectrum in \cite{M20}. 

Motivated by this example, we say that a directed graph $D$ is \emph{spectrally symmetric} if $\Hs(D)$ has spectrum symmetric about $0$. We say that an undirected graph is \emph{spectrally symmetric} if it has at least one \ssym orientation. We focus primarily on oriented graphs because we are interested in the question of density.  

A graph of order $n$ with more than $\frac{n^2}{4}$ edges is not bipartite, hence its adjacency spectrum is not symmetric. We ask for an analogous ``density threshold'' for \ssym orientations. 

\vspace{0.5cm}
\noindent\textbf{Problem 1:} Find the minimal constant $\widehat{\rho}\in [0,1]$ such that each \ssym graph satisfies $\frac{2|E(G)|}{|V(G)|^2}\leq\widehat{\rho}$.
\vspace{0.5cm}

Our main results are a lower bound (Corollary \ref{Cor:lower_density_bound}) and an upper bound (Theorem \ref{Thm:21/22}) showing that \[\frac{31}{36}\leq \widehat{\rho}\le \frac{10}{11}.\]
This bound can be viewed as a spectral analog of Mantel's theorem.

In Section \ref{Sec:A_eigenvalues} we give a necessary condition  (Theorem \ref{Thm:Tr_A^3_min_eval}) for a graph $G$ to admit a \ssym orientation in terms of the eigenvalues of $A(G)$. We apply this condition to show (Theorem \ref{Thm:10n}) that a graph with $|E(G)|>10|V(G)|$ cannot have a spectrally symmetric line graph.

In Section \ref{Sec:Constructions} we show (Theorem \ref{Thm:ind_gluing}) that the 1-sums of certain \ssym graphs are still \ssymend. This allows for the construction of many infinite families of arbitrarily large \ssym graphs.

\subsection{Other angles} 
Hermitian adjacency matrices for digraphs were first studied by Guo and Mohar \cite{GM15} and independently by Liu and Li \cite{LL15}. In both of these papers, the authors used the Hermitian matrix $H_\frac{\pi}{2}$ based on the fourth roots of unity. For that matrix, both papers prove that any oriented graph has spectrum symmetric about 0. So our present question is trivial for $H_\frac{\pi}{2}$.

Subsequently, Mohar \cite{M20} introduced the variant $\Hs$. The theory of $\Hs$ has developed steadily since then. Li and Yu characterized mixed graphs with  $\Hs$-eigenvalues in certain small intervals, \cite{LY21}. Together with Zhou \cite{ZLY23}, they then characterized the mixed graphs with the least $\Hs$-eigenvalue greater than $-\frac{3}{2}$. Van Dam and Wissing studied the closely related generalization to gain graphs defined over the sixth roots of unity, characterizing the examples with rank at most 3 \cite{vDW22b}. Kadyan and Bhattacharjya determined when the eigenvalues of $\Hs$ are integral for certain families of Cayley graphs, \cite{K23,KB22}. The present authors recently determined all mixed graphs with exactly two distinct $\Hs$-eigenvalues \cite{AALMX25}.

The choice of $\theta=\frac{\pi}{3}$ has several advantages over other angles. First, $e^i\theta+e^{-i\theta}=1$, so for mixed graphs $\Hs$ captures algebraically the graph theoretic convention that two arcs in opposite directions should be treated as an undirected edge. Second, $\theta\in \{\frac{2\pi}{3},\frac{\pi}{2},\frac{\pi}{3}\}$ are the only (non-trivial) angles for which the characteristic polynomial of $H_\theta$ is guaranteed to have integer coefficients. This is discussed in \cite{AALMX25}.

The study of \ssym digraphs is nascent. Mohar gave the first example in \cite{M20}. We know of no further results on \ssym orientations for $\Hs$. However, in a different direction, Higuchi, Kubota and Segawa \cite{HKS24} recently found similar examples for which $H_\theta$ is symmetric whenever $\theta\in \pi\mathbb{Q}$. There they also proved that, for irrational $\theta$, only bipartite digraphs have $H_\theta$ symmetric.

\subsection{Preliminary definitions}

A \emph{digraph} $D$ consists of a set $V$ of vertices and a multiset $A$ of ordered pairs of vertices called \emph{arcs}. The \emph{underlying graph} $G$ of $D$ is a multigraph on the same vertex-set with one copy of the edge $\{u,v\}$ for each arc $(u,v)$ or $(v,u)$ in $A(D)$. If $E(G)$ has the same size as $A(D)$ then $D$ is said to be an \emph{oriented} graph. In this case, we will write $E(D)$ instead of $A(D)$ to avoid confusion with the notation for the adjacency matrix. The \emph{in-degree} (respectively, \emph{out-degree}) of a vertex $v$ is the size of the set $\{u: (u,v)\in A\}$ (respectively, $\{u: (v,u)\in A\}$) and is denoted by $d^-(v)$ (respectively, $d^+(v)$). Two vertices $u$ and $v$ in a digraph are \emph{twins} if $d^+(u)=d^+(v)$ and $d^-(u)=d^-(v)$. A digraph is \emph{$d$-regular} if all in-degrees and all out-degrees of all vertices are equal to $d$. A digraph is \emph{bipartite} if its underlying graph is bipartite. A \emph{walk} in an oriented graph is a sequence $v_1,v_2,\dots, v_k$ of vertices such that for each $i\in \{1,\dots,k-1\}$, the pair $v_i$ and $v_{i+1}$ are ends of some arc $a_i$. Note that we allow $a_i=(v_i,v_{i+1})$ and $a_i=(v_{i+1},v_i)$. The walk is \emph{closed} if $v_1=v_k$. A \emph{directed walk} is a walk in which each $a_i=(v_i,v_{i+1})$.
A \emph{triangle} in a digraph is a 3-vertex subgraph with underlying graph $K_3$. A triangle is \emph{directed} if it is regular. A triangle that is not directed is \emph{transitive}. The number of subgraphs of $G$ isomorphic to $K_3$ is denoted by $t(G)$.

\subsection{Small examples}\label{Sec:smallex}
In Sections \ref{Sec:density_bound} and \ref{Sec:Constructions} we will provide methods for constructing new \ssym oriented graphs from given examples, so it is critically important to know that some \ssym orientations of graphs exist. It is easy to see that no examples exist for $n<5$, see the remark after Lemma \ref{Lem:Tr_cubed_tourn}. Here we report on a computer search of the orientations on graphs on five and six vertices which was carried out in Sage \cite{sagemath}. We omit all orientations of  bipartite graphs from the discussion as they are trivially \ssymend.

For $n=5$ only $K_5-E(K_3)$ has \ssym orientations, though there are many orientations with this property.  To discuss all such orientations cogently we need a notion of equivalence of orientations. For simplicity we use a somewhat coarse notion, saying that two orientations of $G$ are \emph{spectrally equivalent} if they have the same $\Hs$-spectrum. There are precisely two spectral equivalence classes of orientations of $K_5-E(K_3)$. A representative for each one is given in Figure \ref{fig_k23}.

\begin{figure}
    \centering
    \includegraphics[scale=1.1]{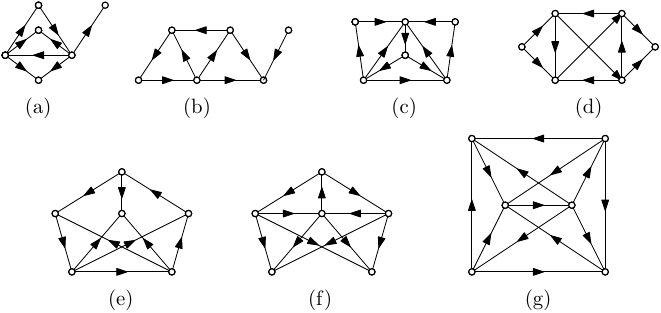}
    \caption{The seven non-bipartite \ssym graphs of order six, each with a \ssym orientation.}
    \label{fig:six_vert_orientations}
\end{figure}

For $n=6$, there are seven \ssym graphs each of which admits different \ssym orientations which belong to at least two spectral equivalence classes. In Figure \ref{fig:six_vert_orientations}, we provide one spectrally symmetric orientation for each of the seven graphs. 

Let us pause to make a few comments on these examples:\begin{itemize}
    \item[(i)] Example (a) is obtained from the 5-vertex example by adding a pendant edge incident with each odd cycle. However, deleting the pendant vertex from (b) does not give a symmetric 5-vertex example because the pendant edge was not incident with each triangle. See Section \ref{Sec:Constructions} for the explanation of this phenomena.
    \item[(ii)] Each of the graphs (a), (b), (e), (f) has a vertex that is contained in all odd cycles. This will be relevant for our construction in Section \ref{Sec:Constructions}.
    \item[(iii)] Examples (d)--(g) each have density greater than $K_{3,3}$. So the density threshold of $\frac{1}{2}$ for undirected graphs with symmetric spectrum is insufficient for oriented graphs. In particular, (e) and (f) are obtained by adding one or two edges to $K_{3,3}$. 
    \item[(iv)] Example (g) is $K_6-E(2K_2)$, it achieves the maximum density among all \ssym oriented graphs that we know of.
\end{itemize}

\section{A density threshold for spectral symmetry} \label{Sec:density_bound}
A graph of order $n$ with more than $\frac{n^2}{4}$ edges is not bipartite, hence its adjacency spectrum is not symmetric. In this section we provide an analogous upper bound on the density of an oriented graph with symmetric $\Hs$-spectrum. 

The typical notion of density for a graph (or oriented graph) of order $n$ is the ratio $|E(G)|/|E(K_n)|$. For our purpose, it is more natural to work with the \emph{limit density} defined by \[\widehat{\rho}(G):=\frac{2|E(G)|}{|V(G)|^2},\] so called because it approaches the aforementioned ratio as $|V(G)|\rightarrow\infty$. We prefer this asymptotic notion of density because we can construct arbitrarily large \ssym graphs with limit density exactly $\widehat{\rho}(G)$ from a given example $G$. As an added bonus, the limit density of $K_{n,n}$ is exactly $\frac{1}{2}$, which seems appropriate.

\subsection{The upper bound}

First we need a simple lemma. 

\begin{lemma} \label{Lem:Tr_cubed_tourn}
Suppose that $D$ is a \ssym orientation. Let $t_0$ be the number of transitive triangles in $D$ and $t_1$ be the number of directed triangles. Then $t_0 = 2t_1$. 
\end{lemma}

\begin{proof}
Let $H=\Hs(D)$. Note that $\tr(H^3)$ counts triples $(v_1,v_2,v_3)$ that form a
triangle. In this counting each transitive triangle is counted three
times with weight $\frac{1+i\sqrt{3}}{2}$ and three times with weight
$\frac{1-i\sqrt{3}}{2}$, giving a contribution of 3 altogether. On the other
hand, each directed triangle is counted six times with weight $-1$ each
time. Since $H$ has symmetric spectrum, we have $\tr(H^3)=0$, which is
equivalent to $3t_0 - 6t_1 = 0$. Therefore, we must have $t_0=2t_1$.    
\end{proof}

\noindent\textbf{Remark:} As an immediate corollary, a graph that admits a \ssym orientation must have a number of triangles divisible by three. From this it is easy to see that any non-bipartite \ssym graph has at least five vertices.

\begin{theorem}[Upper bound on $\widehat{\rho}$]\label{Thm:21/22}
    If $D$ is a \ssym orientated graph, then $\widehat{\rho}(D)\le\frac{10}{11}.$
\end{theorem}

\begin{figure}[ht]
    \centering
    \includegraphics[width=0.5\linewidth]{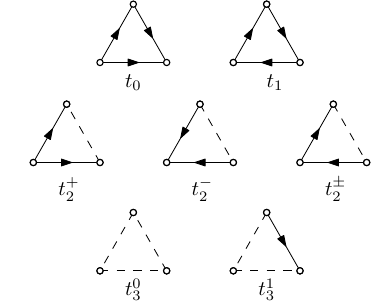}
    \caption{The induced subgraphs of order 3 counted in the proof of Theorem \ref{Thm:21/22}. Dashed lines indicate edges of $R$.}
    \label{fig:triangle_types}
\end{figure} 

\begin{proof}
Let $n$ be an arbitrarily large number divisible by $|V(D)|$. By Theorem \ref{Thm:Density_Lemma} the $\frac{n}{|V(D)|}$-fold blowup of $D$ is \ssym with the same limit density. It suffices to show that this oriented graph--call it $B$--of order $n$ has limit density at most $\frac{10}{11}$. Let $R$ denote the complement of the underlying graph of $B$ and let $d_v^R$ denote the degree of $v$ in $R$. Let $r:=|E(R)|$.

Let $t_0$ and $t_1$ denote the number of transitive and directed triangles in $B$, respectively. Moreover, let $t^{+}_2, t^-_2, t^{\pm}_2, t^0_3, t^1_3$ denote the number of oriented induced subgraphs of order 3 in $B$ of the types indicated in Figure \ref{fig:triangle_types}.

Since $B$ has symmetric spectrum, Lemma \ref{Lem:Tr_cubed_tourn} implies $t_1=\frac{1}{2}t_0$. The total number of induced subgraphs of order 3 in $B$ is $ \binom{n}{3}$ and so the following holds: \begin{equation} \label{eq_2.4_first}
    t_0+  t_1+t^{+}_2+ t^-_2+t^{\pm}_2+t^0_3+ t^1_3=  \binom{n}{3}= \frac{n^3}{6}+O(n^2).
\end{equation}

Also, one can easily verify the following: $$t_0=  \sum_v  \binom{d^+_v}{2}- t^+_2= \sum_v \binom{d^-_v}{2}- t^-_2.$$

By considering  half of the sum of  the left and  right sides of the previous equality, we have  $$t_0=\frac{1}{4}\sum_v [(d^+_v)^2 + (d^-_v)^2]-\frac{1}{2} (t^+_2 +  t^-_2)-O(n^2).$$

Now, since $t_0+t_1=\frac{3}{2}t_0$, using Equation (\ref{eq_2.4_first}), one can see that 

\begin{equation}\label{eq_2.4_second}
n^3= \frac{9}{4}\sum_v [(d^+_v)^2 + (d^-_v)^2]+\frac{3}{2} (t^+_2 + t^- _2) +6 (t^{\pm}_2+ t^0 _3+t^1 _3)+O(n^2).
\end{equation}

In the sequel, we shall use the following set:  $$\mathcal{A} =\{ \, (e, u) \,\, |\, e=xy \in E(R), u\in V(G), \,\,  u\neq x, y \,\,\}.$$ Now, we would like to count the size of $\mathcal{A}$ in two different ways. For any edge $e$ of $R$, there are $n-2$ vertices which are not incident with $e$. Thus we have  $|\mathcal{A}|=r(n-2)$. On the other hand, by considering the induced subgraphs of order $3$ and fixing a vertex $u$,  one can see that $$|\mathcal{A}|=t^{+}_2+ t^-_2+t^{\pm}_2+ 2t^1_3 + 3 t^0_3.$$

Now, using Equation (\ref{eq_2.4_second}), we obtain that 

\begin{equation} \label{eq_2.4_third}
    n^3(1+o(1))=\frac{9}{4}\sum_v [(d^+_v)^2 + (d^-_v)^2]+\frac{3}{2}rn+\frac{9}{2}t^{\pm}_2 + 3t^1_3 + \frac{3}{2} t^0_3.
\end{equation}

By the Cauchy–Schwarz Inequality and equation $\sum_v d^R_v=2r$, we have

\begin{align*}
3t^1_3 + \frac{3}{2} t^0_3&\geq \frac{1}{2}(t^1_3+3t^0_3)\\&=\frac{1}{2} \sum_v  \binom{d^R_v}{2}\\&= \frac{1}{4} \sum (d^R_v)^2-\frac{1}{2}r\\&\geq \frac{1}{4n}(2r)^2-\frac{1}{2}r\\&=\frac{r^2}{n}-\frac{1}{2}r.    
\end{align*}

Together with this inequality, equation (\ref{eq_2.4_third}) yields $$ n^3(1+o(1))\geq \frac{9}{4}\sum_v [(d^+_v)^2 + (d^-_v)^2]+\frac{3}{2}rn+ \frac{r^2}{n}-\frac{1}{2}r.$$

Also, we have  $\sum_v (d^+_v + d^-_v)=2|E(G)|=2 (\binom{n}{2}-r)= n^2-2r-n.$

By the Cauchy–Schwarz Inequality and the previous approximation, the following holds: \begin{align*}
\sum_v [(d^+_v)^2 + (d^-_v)^2]&\geq \frac{1}{2n} (n^2-2r-n)^2\\&=\frac{1}{2}n^3-2rn+\frac{2r^2}{n}-n^2+2r+\frac{n}{2} \\&\geq\frac{1}{2}n^3(1+o(1))-2nr
+\frac{2r^2}{n}.\end{align*}

Let $\alpha=\frac{r}{n^2}$. Since $n$ is arbitrarily large, we conclude that $$n^3\geq \frac{9}{8}n^3-\frac{9}{2}rn+\frac{9}{2}\Big(\frac{r^2}{n}\Big)+\frac{3}{2}rn+\frac{r^2}{n},$$ which implies that $\frac{1}{8} n^3-3rn +\frac{11}{2}\frac{r^2}{n}\leq0,$ or equivalently we have the inequality $ \frac{1}{8}-3\alpha+\frac{11}{2}\alpha^2\leq 0. $
Therefore, we find that $$\alpha \geq \frac{24-\sqrt {24^2-4(44)}}{88}=\frac {1}{22}.$$ So \[\widehat{\rho}(D)=\widehat{\rho}(B)=\frac{n(n-1)-2r}{n^2}=1-\frac{1}{n}-2\alpha\leq \frac{10}{11},\] and the proof is complete.
\end{proof}

In particular, Theorem \ref{Thm:21/22} shows that no tournament of order $n>11$ is \ssymend. To conclude this section, we show that, in fact, no tournament is \ssymend.

\begin{theorem} \label{Thm:no_sym_tourn}
No tournament of order $n\ge 3$ is \ssymend.
\end{theorem}

\begin{proof}
Let $T$ be a tournament of order $n\ge3$ and suppose, for a
contradiction, that $\Hs(T)$ has symmetric spectrum.

Let us first observe that $T$ has $\binom{n}{3}$ triangles. Since every
transitive triangle has precisely one vertex with two outgoing edges, it
is clear that the number of transitive triangles is equal to:
\[t_0 = \sum_{i=1}^n \binom{d^+_i}{2}.\] 
By Lemma \ref{Lem:Tr_cubed_tourn} we conclude that

\begin{equation} \label{eq:Tourn_proof_2nd}
     \sum_{i=1}^n \binom{d^+_i}{2} = \frac{2}{3}\binom{n}{3}.
\end{equation}
After using the fact that $\sum_{i=1}^n d^+_i = \binom{n}{2}$, this
equation simplifies to

\begin{equation}\label{eq:Tourn_proof_3rd}\frac{4}{3}\binom{n}{3} + \binom{n}{2} = \sum_{i=1}^n (d^+_i)^2 \ge
\frac{1}{n}{\binom{n}{2}}^2,\end{equation}
where the last inequality uses the Cauchy-Schwartz inequality. Note also
that the last inequality is strict unless $T$ is regular, i.e.
$d^+_i=\frac{n-1}{2}$ for every $i$.

The last inequality simplifies to the condition that $n\le11$. Equality
(\ref{eq:Tourn_proof_2nd}) implies that one of $n$, $n-1$ or $n-2$ must be
divisible by $9$. This leaves only the cases $n=9,10,11$. For these
values, the left-hand side in (\ref{eq:Tourn_proof_3rd}) is equal to $148$, $205$, and $275$, respectively. It is elementary to check that the only way for 
$n$ integers to sum to $36$, $45$, $55$ (respectively), while their squares sum to $148$, $205$, and $275$ (respectively) are: $\{5^{(2)},4^{(5)},3^{(2)}\}, \{5^{(5)},4^{(5)}\}, \{5^{(11)}\}$, respectively.

The case for $n=11$ is that of a regular tournament, and in fact, a short argument shows that no regular tournament has symmetric spectrum: A regular tournament has odd order, so if it has symmetric spectrum, then $0$ is an eigenvalue of $\Hs(T)$. Let $v$ be the eigenvector of $0$. Note that $\Hs(T)=\frac{1}{2}(J-I+iS)$. Since $T$ is regular, $S$ and $T$ are simultaneously diagonalizeable and $\frac{n-1}{2}$ is the eigenvalue for the all ones eigenvector. All other eigenvectors, including $v$, are orthogonal to $J$ so we have \[Iv=iSv=i\sqrt{3}S^\prime v\] with $S^\prime$ a skew symmetric matrix with integer entries. So the eigenvalues of $S^\prime$ must be algebraic integers, but the above equation implies $S^\prime v=-\frac{i}{\sqrt{3}}v$, a contradiction.

For the remaining cases of $n=9,10$ with the specified degrees we could not find a simple proof, so we checked computationally that none exist. We used the catalogue of directed graphs provided by Brendan McKay \cite{BMweb} and computed the $\Hs$-eigenvalues of the relevant tournaments. There are $3,310$ tournaments of order $9$ to consider and $13,333$ of order $10$, each was found to have asymmetric spectrum. 
\end{proof}

In the proof above, it is tempting to consider the $n=9,10$ cases without computer aid. For instance, one could hope that $\tr(H^5)=0$ provides additional constraints on the orientation to rule out these cases. However, there is a tournament on $9$ vertices with the curious property that each of $H^3,H^5,H^7$ has trace $0$ and only $\tr(H^9)\neq0$ prevents this tournament from having symmetric spectrum.

\subsection{The lower bound}

To bound $\hat{\rho}$ from below we wish to construct a very dense spectrally symmetric orientation. As a first candidate, the digraph in in Figure \ref{fig:six_vert_orientations} (g) gives a a lower bound of $\frac{13}{18}$ on $\hat{\rho}$.

In a previous version of this paper posted to the arxiv, we suggested the improvement of this bound as an open problem. Seamons found such an improvement, using computer search to find a spectrally symmetric orientation of $K_{11}-3K_2$ in \cite{JS26}. This dramatically increases the lower bound from $\frac{13}{18}$ to $\frac{104}{121}$, which is quite close to the $\frac{10}{11}$ upper bound of the previous section. Emboldened by this success, we have found an even denser spectrally symmetric orientation of $K_{12}-4K_2$, pictured in Figure \ref{fig:K12-4K2}.

\begin{figure}
    \centering
    
\begin{tikzpicture}[
    x=1cm,y=1cm,
    orient/.style={
        draw=black,
        line width=0.55pt,
        postaction={decorate},
        decoration={
            markings,
            mark=at position 0.56 with
                {\arrow{Latex[length=2.8mm,width=2.2mm]}}
        }
    },
    wtone/.style={
        circle,
        draw=black,
        fill=white,
        line width=0.65pt,
        inner sep=0pt,
        minimum size=4pt
    },
    wttwo/.style={
        circle,
        draw=black,
        fill=white,
        line width=0.65pt,
        inner sep=0pt,
        minimum size=8pt
    },
    vlabel/.style={font=\small}
]

\def\R{2.65}

\coordinate (A) at ( 90:\R);
\coordinate (p) at ( 45:\R);
\coordinate (B) at (  0:\R);
\coordinate (q) at (-45:\R);
\coordinate (C) at (-90:\R);
\coordinate (r) at (-135:\R);
\coordinate (D) at (180:\R);
\coordinate (s) at (135:\R);

\draw[orient] (A) -- (B);
\draw[orient] (A) -- (C);
\draw[orient] (A) -- (q);

\draw[orient] (B) -- (C);
\draw[orient] (B) -- (r);
\draw[orient] (B) -- (s);

\draw[orient] (C) -- (p);
\draw[orient] (C) -- (r);
\draw[orient] (C) -- (s);

\draw[orient] (D) -- (A);
\draw[orient] (D) -- (B);
\draw[orient] (D) -- (C);
\draw[orient] (D) -- (p);

\draw[orient] (p) -- (A);
\draw[orient] (p) -- (B);
\draw[orient] (p) -- (q);
\draw[orient] (p) -- (r);
\draw[orient] (p) -- (s);

\draw[orient] (q) -- (B);
\draw[orient] (q) -- (C);
\draw[orient] (q) -- (D);

\draw[orient] (r) -- (A);
\draw[orient] (r) -- (D);
\draw[orient] (r) -- (q);
\draw[orient] (r) -- (s);

\draw[orient] (s) -- (A);
\draw[orient] (s) -- (D);
\draw[orient] (s) -- (q);

\node[wttwo,label={[vlabel,label distance=2.5pt]90:$A$}] at (A) {};
\node[wtone,label={[vlabel,label distance=2.5pt]45:}] at (p) {};
\node[wttwo,label={[vlabel,label distance=2.5pt]0:$B$}] at (B) {};
\node[wtone,label={[vlabel,label distance=2.5pt]-45:}] at (q) {};
\node[wttwo,label={[vlabel,label distance=2.5pt]-90:$C$}] at (C) {};
\node[wtone,label={[vlabel,label distance=2.5pt]-135:}] at (r) {};
\node[wttwo,label={[vlabel,label distance=2.5pt]180:$D$}] at (D) {};
\node[wtone,label={[vlabel,label distance=2.5pt]135:}] at (s) {};

\end{tikzpicture}

    \caption{Replacing the large vertices $A,B,C,D$ with two twin vertices each yields a spectrally symmetric orientation of $K_{12}-4K_{2}$.}
    \label{fig:K12-4K2}
\end{figure}

We identified $K_{12}-4K_2$ as a likely candidate for a dense spectrally symmetric orientation, and prescribed the condition that each pair of non-adjacent vertices should be twins. This ensures $0$ is an eigenvalue of multiplicity 4 in the Hermitian adjacency matrix and means we can search over weighted tournaments graphs of order $8$ with four vertices of weight 2. This smaller search space contained the pictured example.

We also note that any example of a digraph of density $\hat{\rho}$ extends to an infinite family by taking blowups. 

\begin{theorem}\label{Thm:Density_Lemma}
    Suppose $D$ is a \ssym digraph with limit density $\widehat{\rho}=\widehat{\rho}(D)$. Let $J_\ell$ be the $\ell\times \ell$ matrix of ones and let $D_\ell$ be the digraph whose Hermitian adjacency matrix is the Kronecker product $\Hs(D) \otimes J_\ell$. Then $D_\ell$ is a \ssym orientation with limit density $\widehat{\rho}$.
\end{theorem}

\begin{proof}
     Note that $D_\ell$ is a blowup of $D$ obtained by replacing each vertex $v$ with an independent set $I_v$ of size $\ell$ and replacing each arc $(u,v)$ with a copy of $K_{\ell,\ell}$ between $I_u$ and $I_v$ so that all arcs are oriented from $u$ to $v$. The eigenvalues of $\Hs(D)\otimes J_\ell$ are the products of the eigenvalues of $\Hs(D)$ and those of $J_\ell$. Therefore, since $D$ has symmetric spectrum, so does $\Hs(D)\otimes J_\ell$.

     Let $n$ and $m$ denote the number of vertices and arcs in $D$. The limit density of $D_\ell$ is 
    \[\frac{2|E(D_\ell)|}{|V(D_\ell)|^2}=\frac{2m\ell^2}{(n\ell)^2}=\frac{2m}{n^2}=\widehat\rho(D).\qedhere\]
\end{proof}

Note that the matrix $J_\ell$ could be replaced with any symmetric $(0,1)$-matrix. The product would still give rise to a symmetric digraph. For us, $J_\ell$ is ideal since it gives the largest density.

\begin{corollary}[Lower bound on $\widehat{\rho}$] \label{Cor:lower_density_bound}
    There are infinitely many \ssym oriented graphs with limit density $\frac{31}{36}$.
\end{corollary}

\begin{proof}
    Take blowups of the orientation of $K_{12}-4K_2$ from Figure \ref{fig:K12-4K2}.
\end{proof}

\section{A condition on the adjacency eigenvalues} \label{Sec:A_eigenvalues}

Here we give another necessary condition on graphs that admit \ssym orientations. Remarkably, it can be stated entirely in terms of the adjacency spectrum of the graph. We first need a lemma.

\begin{lemma} \label{Lem:Trace_PSD_ineqality} 
    Suppose $M$ is a symmetric matrix and $X$ is a positive semi-definite matrix. Then $\tr(MX)\leq\lambda_{max}(M)\tr(X)$.
\end{lemma}

\begin{proof}
    $X$ decomposes as $X=\sum_{k=1}^r\mu_k\alpha_k\alpha_k^T$ with $\mu_k\geq0$ the eigenvalues of $X$. So 
    \begin{align*} \tr(MX)&=\sum_{k=1}^r\tr(M\mu_k\alpha_k\alpha_k^T)\\&=\sum_{k=1}^r\mu_k\tr(\alpha_k^TM\alpha_k)\\&\leq\sum_{k=1}^r\mu_k\lambda_{max}(M)\\&=\lambda_{max}(M)\tr(X). \qedhere
    \end{align*} 
\end{proof}

Now, note that 
\begin{equation}
    \Hs(D)=H=\frac{1}{2}(A+iS) \label{Eq:HAS}
\end{equation} 

\noindent where $A$ is the adjacency matrix of the graph underlying $D$ and $S$ is a skew-symmetric matrix with with the same pattern as $A$ and off diagonal entries in $\{0,\pm\sqrt{3}\}.$

\begin{theorem} \label{Thm:Tr_A^3_min_eval}
Let $A=A(G)$ be the adjacency matrix of an undirected graph $G$ with $m$ edges. If $G$ has a \ssym orientation, then \[\tr(A^3)\le-18m\lambda_{min}(A).\]
\end{theorem}

\begin{proof}
If $H=\Hs(D)$ has symmetric spectrum then $\tr(H^3)=0$. We use Equation (\ref{Eq:HAS}) and expand $H^3$ as a sum of 4 terms. By the cyclic property of the trace we have
\begin{align*}
0=\tr(H^3)&=\frac{1}{8}\tr\Big((A+iS)^3\Big)\\&=\frac{1}{8}\tr(A^3)-\frac{3}{8}\tr(AS^2)+\frac{3i}{8}\tr(A^2S)-\frac{i}{8}\tr(S^3)\\&=\frac{1}{8}\tr(A^3)-\frac{3}{8}\tr(AS^2).    
\end{align*}

Here the final equation follows from the fact that $S^3$ is skew symmetric and $\tr(A^2S)=\tr((A^2S)^T)=\tr(S^T(A^T)^2)=-\tr(A^2S)=0.$ 

Let $d_v$ be the degree of vertex $v$ of $G$. Note that $A$ is symmetric, $-S^2$ is positive semi-definite, and $S^2_{v,v}=-3d_v$ for $v$ in $V(G).$ So Lemma \ref{Lem:Trace_PSD_ineqality} gives  

\begin{align*}    
\tr(A^3)&=3\tr(AS^2)\\&=3\tr((-A)(-S^2))\\&\le-3\lambda_{min}(A)\tr(-S^2)\\&=-18m\lambda_{min}(A). \qedhere
\end{align*}

\end{proof}

Theorem \ref{Thm:Tr_A^3_min_eval} gives a nice corollary for line graphs, since the (adjacency) eigenvalues of a line graph are bounded below by $-2$.

\begin{corollary}\label{Cor:t<6e}
If $L$ is a \ssym line graph then
\[t(L)\le6|E(L)|.\]  
\end{corollary}

This shows that the line graphs of reasonably dense graphs cannot have spectrally symmetric orientations.

\begin{theorem}\label{Thm:10n}
    If the line graph $L=L(G)$ is \ssymend, then $|E(G)|\le10|V(G)|$.\end{theorem}

\begin{proof}

Let $G$ be a graph with $n$ vertices. Suppose $L$ is \ssymend. Then using Corollary \ref{Cor:t<6e} and counting triangles and edges of $L$ from the degrees of $G$ we have

\[t(G)+\sum_{i=1}^n\binom{d_i}{3}\le6\sum_{i=1}^n\binom{d_i}{2}.\]
After simplification and ignoring the $t(G)$ term we have
\begin{align}\sum\limits_{i=1}^nd_i^3-21d_i^2+20d_i<0\hspace{10pt} \label{eq_cubic}
\end{align}
So it suffices to show that  (\ref{eq_cubic}) fails when $|E(G)|> 10n$. Let $\alpha$ denote the average degree of $G$. Suppose $|E(G)|> 10n$ and thus $\alpha>20$. Then

\begin{align}\sum\limits_{i=1}^n(d_i^3-21d_i^2+20d_i)=&\sum_{i=1}^n((d_i-7)^3-127d_i+343)\nonumber\\ 
\geq&\sum_{i=1}^n((\alpha-7)^3-127\alpha+343)\label{eq_alpha_cubic}\\
=&\sum\limits_{i=1}^n(\alpha^3-21\alpha^2+20\alpha)\nonumber\\
=&~n(\alpha-20)(\alpha-1)\alpha \nonumber\\
>&~0.\nonumber\end{align}

The inequality on Line (\ref{eq_alpha_cubic}) follows because, for $\alpha>14$, $\sum_{i=1}^n(d_i-7)^3$ is minimized when all $d_i$ are equal.
\end{proof}

 \noindent \textbf{Remark:} In the above proof we did not consider the contribution of the $t(G)$ term, so it is possible the constant $10$ can be improved a bit. However $t(L(K_{19}))<6|E(L(K_{19})|$, so we would need other methods to push the constant much further.

\section{Construction of infinite families} \label{Sec:Constructions}
Above we have shown certain limits to the class of \ssym graphs. In this section we provide a construction showing that these graphs are, nevertheless, plentiful. Another goal of this section is to explain the phenomenon mentioned in remark $(i)$ in Section \ref{Sec:smallex}.

We say that a vertex $v$ in a graph $G$ is \emph{odd-dominating} if each odd cycle of $G$ contains $v$. A vertex in a digraph is \emph{odd-dominating} if it is odd-dominating in the underlying graph. Note that the digraphs from Figure \ref{fig_k23} as well as digraphs (a), (b), (e), (f) from Figure \ref{fig:six_vert_orientations} have odd-dominating vertices. 

The main result of this section is that spectral symmetry is preserved when two spectrally symmetric graphs are glued together at a pair of odd dominating vertices. We need a bit more setup to make this statement formal. 

If $v\in V(D_1)$ and $u\in V(D_2)$ are vertices of vertex disjoint digraphs $D_1$ and $D_2$, then the \emph{1-sum} of $D_1$ and $D_2$ at $u$ and $v$ is the graph obtained from $D_1\cup D_2$ by identifying $u$ and $v$. 

Let $\{D_i\}_{i\in \{1,\dots,\ell\}}$ be a sequence of vertex disjoint digraphs, each of which has an odd-dominating vertex. We say that a vertex from any of the $D_i$ is \emph{locally odd-dominating} if it is odd-dominating in some $D_i$. We let $B_1=D_1$ and for $j>1$ we construct $B_j$ as a 1-sum of $B_{j-1}$ with $D_j$ at a pair of locally odd-dominating vertices. Now we can state the main theorem of this section.

\begin{theorem} \label{Thm:Gluing_Construction_1}
If each  $\{D_i\}_{i\in \{1,\dots,\ell\}}$ is \ssymend, then each $B_j$ $(1\leq j\leq \ell)$ is \ssym as well. 
\end{theorem}

In order to verify that our constructions produces \ssym digraphs, we make use of the following folklore result about spectrally symmetric matrices.
\begin{lemma} \label{Lem:Sym_spec_Char}
    Let $H$ be a Hermitian matrix. We have $\tr(H^r)=0$ for all positive odd integers $r$ if and only if $H$ has symmetric spectrum. \qed
\end{lemma}

Let $H=\Hs(D)$. For a walk $W=w_0w_1\cdots w_r$ of length $r$ in $D$, the \emph{value} of $W$ is \[|W|:=\prod_{i=1}^rH_{w_{i-1}w_i}.\] The value of a set of walks $\W$ is $|\W|:=\sum_{W\in \W}|w|.$ Let $\W_r(D)$ denote the set of all closed walks of length $r$ in $D$. For a vertex $v$ in $D$ we define $\W_r[D,v]$ (and $\W_r(D,v)$) to be the walks in $\W_r(D)$ which visit $v$ at least (exactly) once. We will sometimes suppress $D$ from this notation.

Note that $\tr(H^r)=|\W_r(D)|$. So, to prove that $D$ is \ssym it suffices to show that $|\W_r(D)|=0$ for all odd $r$, and then apply Lemma \ref{Lem:Sym_spec_Char}. Our main technical tool in this endeavor is a special partition of the odd walks through a given vertex $v$.

Let $\S\subseteq \W_r[D,v]$. The \emph{first-odd-subwalk partition (fos-partition) of} $\S$ is defined in three steps. 

\begin{enumerate}
    \item We partition $\S$ into sets $\S_{x,y}$ as follows: For each closed walk $W=w_0w_1\dots w_{r}$ $(w_0=w_r)$ find the lexicographically least pair of indices $(x,y)$ such that $x$ is odd and $w_y=w_{y+x~\bmod{r}}=v$. All walks in $\S$ are admitted into some part because $r$ is odd. By the lexicographical condition $w_{z}\not=v$ for $y<z<y+x$, and thus the subwalk $W^\prime=w_yw_{y+1}\ldots w_{y+x}$ of $W$ is an element of $\W_x(D,v)$. Note that all indices are mod $r$.

    \item If $D=B_j$ is one of the iteratively constructed graphs above, and if $v$ is one of the vertices at which the constituent pieces $D_i$ were glued together, then we further decompose each $\S_{x,y}$ into sets $\S_{x,y,i}$ according to which of the subgraphs $D_i$ contains the first odd subwalk $W'$. If $D$ is not one of the iteratively constructed graphs, simply skip this step.

    \item Let $\circ$ denote the concatenation of walks. If $x+y\leq r$ then each element of $\S_{x,y,i}$ is of the form $W_b\circ W^\prime\circ W_e$ where $W^\prime\in \W_{x}(D,v)$. Here we further decompose each $\S_{x,y,i}$ into sets $\P_{x,y,i}(W_b,W_e)$ with beginning subwalk $W_b$ and ending subwalk $W_e$ and any element of $\W_{x}(D,v)$ in between them.
    
    If $x+y>r$ then $W$ is of the form $W_2^\prime\circ W_f\circ W_1^\prime$ where $W_1^\prime\circ W_2^\prime\in \W_x(D,v).$ These differ from the previous case only in notation, so we will not continue to discuss them explicitly. We also suppress the indices of $\P$ for readability, so a fos-partition of $\S$ is just $\{\P(W_b,W_e)\}.$
\end{enumerate}

Each part $\P(W_b,W_e)$ consists of a fixed beginning subwalk and a fixed ending subwalk. The variable part in between is a set of subwalks of $\S$ of fixed odd length and through a fixed $D_i$ subgraph if $D=B_j$. The reason for this lengthy definition is that it will be easy to show inductively that the value of this variable set of walks is $0$, hence each $|\P(W_b,W_e)|=0$, and so each $|\S_{x,y,i}|=0$ as well.

\begin{lemma}\label{Lem:odd_walk_decomp}
Let $D$ be a digraph and $v\in V(D)$ be a vertex. Then $|\W_r[D,v]|=0$ for all odd $r$ if and only if $|\W_r(D,v)|=0$ for all odd $r$.   
\end{lemma}

\begin{proof}
First suppose $|\W_r(D,v)|=0$ for all odd $r$. We construct the fos-partition of $\W_r[D,v]$.

Note that each element of $\P(W_b,W_e)$ has a fixed contribution of $|W_b|\cdot|W_e|$ and a variable contribution of $|W^\prime|$ where $W^\prime$ varies over elements of $\W_x(v)$ which begin at $v$. The value of all such walks is $\frac{1}{x}|\W_x(v)|$, so \[|\P(W_b,W_e)|=\frac{1}{x}|\W_x(v)|\cdot|W_b|\cdot|W_e|.\] 
Since $|\W_{x}(v)|=0$, we are done.

The other direction follows by induction on $r$. The base case is $r=1$. Clearly, $\W_1(D,v)$ is empty so has value zero. Now, for the induction suppose $|\W_s(D,v)|=0$ for all odd $s<r$ and $|\W_r[D,v]|=0$. Apply the fos-partition on the set $\W_r[D,v]\setminus\W_r(D,v)$. Because $v$ occurs at least twice, $x<r$ for each $\S_{x,y}$, so each of the constituent parts $\P(W_b,W_e)$ has value $0$ by our induction hypothesis. Thus \[|\W_r(D,v)|=|\W_r[D,v]|-|\W_r[D,v]\setminus\W_r(D,v)|=0.\qedhere\]  
\end{proof}

Now we state and prove an inductive version of Theorem \ref{Thm:Gluing_Construction_1}.

\begin{theorem}\label{Thm:ind_gluing}
Let $D_1,\dots,D_\ell$ be \ssym and let $B_j$ $(1\leq j\leq \ell)$ be constructed as before the statement of Theorem \ref{Thm:Gluing_Construction_1}. Then the following holds:
\begin{enumerate}
    \item[(i)] $B_\ell$ is \ssymend.
    \item[(ii)] For any vertex $u$ of $B_\ell$ which is locally odd-dominating, and for each odd integer $r$, $|\W_r(D,u)|=0$.
\end{enumerate}
\end{theorem}

\begin{proof}
We proceed by induction on $\ell$. The first claim is our main goal, the second is an aid in our induction. When $\ell=1$, claim $(i)$ is immediate and claim $(ii)$ follows from Lemma \ref{Lem:odd_walk_decomp} since $\W_r=\W_r[B_1,u]$. 

Now for $\ell>1$ suppose $(i)$ and $(ii)$ hold for each index $j<\ell$. We first prove $(i)$. Let  $v\in V(B_\ell)$ be the vertex at which $D_\ell$ and $B_{\ell-1}$ were glued. For each odd $r$, the $r$-walks of $B_\ell$ are of three types, $\W_r(B_\ell)=\T_1\cup \T_2\cup \T_3$. Walks in $\T_1$ never leave $D_\ell$, walks in $\T_2$ never leave $B_{\ell-1}$ and walks in $\T_3$ cross between the parts at least once at the vertex $v$. Then $|\T_1|=|\T_2|=0$ because $D_\ell$ and $B_{\ell-1}$ are \ssym by $(i)$ inductively. It remains to show that $|\T_3|=0$. 

We employ the fos-partition of $\T_3$ at $v$. We find that the value of $\T_3$ is a sum of multiples of $|\W_{r^\prime}(B_{\ell-1,},v)|$ and $|\W_{r^\prime}(D_{\ell},v)|$ for various odd $r^\prime<r$, and these are $0$ by $(ii)$ inductively. Claim (i) now follows from Lemma~\ref{Lem:Sym_spec_Char}.

Now we will prove $(ii)$. Let $u$ be some locally odd dominating vertex of $B_\ell$, and $r$ an odd integer. We may assume there is a unique index $j\leq \ell$ such that $u\in V(D_j)$ since otherwise we simply conclude with the argument from the previous paragraph with $\W_r(B_\ell,u)$ in place of $\T_3$.

Let $\overline{\W_r[B_\ell,u]}$ be the complement of $\W_r[B_\ell,u]$ in $\W_r(B_\ell)$, i.e.~the closed walks of length $r$ in $B_\ell$ which miss $u$. We next show that $|\overline{\W_r[B_\ell,u]}|=0$. Note that each walk in $\overline{\W_r[B_\ell,u]}$ traverses the edges of $D_j$ an even number of times, since otherwise the walk would have an odd closed subwalk contained in $D_j$ but missing $u$. So, the value of $|\overline{\W_r[B_\ell,u]}|$ will be a multiple of a sum $\sum_{r'<r}|\W_{r'}[B',u']|$ with $r'$ odd and $B^\prime$ obtained from $B_\ell$ by contracting all vertices of $D_j$ into a single vertex $u^\prime$. Since $B'$ is obtained from a subset of the $\{D_i\}$ of size $\ell-1$ the terms of this sum are $0$ by $(ii)$ inductively and by Lemma \ref{Lem:odd_walk_decomp}. Hence $|\overline{\W_r[B_\ell,u]}|=0$.

Finally, since $(i)$ holds, we have  $|\W_r(B_\ell)|=0$. So~$|\W_r[B_\ell,u]|=|\W_r(B_\ell)|-|\overline{\W_r[B_\ell,u]}|=0-0$. Now $(ii)$ follows from Lemma \ref{Lem:odd_walk_decomp}.
\end{proof}

\noindent\textbf{Remark:} We did not need simple digraphs $D_i$ for the construction to work. Nor did we need to use the matrix $\Hs(D)$. Any Hermitian adjacency matrix $H_\theta(D)$ would work, provided we can find small examples to build from. One such example for each $\theta \in \pi \mathbb{Q}$ was given by Higuchi Kubota and Segawa in \cite{HKS24}. Generalizing the example from Figure \ref{fig_k23}, they produced for each such $\theta$ a non-bipartite oriented graph for which $H_\theta$ has symmetric spectrum. It is easy to see that their examples, like $K_5-E(K_3)$, have odd-dominating vertices, so we can use those examples as seeds to grow a variety of digraphs with symmetric $H_\theta$-spectrum.

\section*{Statement on the use of AI}
The spectrally symmetric orientation of $K_{12}-4K_2$ from Figure \ref{fig:K12-4K2} was found by ChatGPT Sol 5.6 (the agent) during the course of a conversation with one of the authors. The agent identified $K_{12}-4K_2$ as a viable candidate for improving the density bound, and the author suggested copying the structure of the previously known orientation of $K_{6}-2K_2$. The agent recognized the two pairs of twin vertices in that example and correctly extrapolated to find the new example. The authors take full responsibility for the correctness of this example, and an auxiliary file verifying its spectral symmetry using sage is available upon request. With this notable exception, no AI tools were used in course of this research or in the preparation of this manuscript.

\bibliographystyle{abbrv}
\bibliography{HermMat}

\end{document}